\documentclass[11pt]{article}
\RequirePackage[OT1]{fontenc}
\RequirePackage{amsthm,amsmath}
\RequirePackage[numbers]{natbib}

\newcommand{\rr}{{\Bbb R}}

\newcommand{\one}{{\hbox{1{\kern -0.35em}1}}}

\newcommand{\M}{{\cal M}}

\newcommand{\beq}[1]{\begin{eqnarray} \label{#1}}
\newcommand{\eeq}{\end{eqnarray}}
\newcommand{\bed}{\begin{displaymath}}
\newcommand{\eed}{\end{displaymath}}
\newcommand{\bea}{\bed\begin{array}{rl}}
\newcommand{\eea}{\end{array}\eed}
\newcommand{\disp}{\displaystyle}
\newcommand{\al}{\alpha}

\newtheorem{thm}{Theorem}[section]

\newtheorem{lem}[thm]{Lemma}
\newtheorem{rem}[thm]{Remark}

\newtheorem{defn}[thm]{Definition}

\newcommand{\thmref}[1]{Theorem~{\rm \ref{#1}}}
\newcommand{\lemref}[1]{Lemma~{\rm \ref{#1}}}

\def\openbox{$\sqcup\llap{$\sqcap$}$}
\def\endproof{\unskip \enskip
    \null \nobreak \hfill \openbox \par}

\newcommand{\ad}{&\!\!\!\disp}
\newcommand{\aad}{&\disp}
\newcommand{\barray}{\begin{array}{ll}}
\newcommand{\earray}{\end{array}}

\usepackage{makeidx}  
\usepackage{amsfonts}
\makeindex            
\usepackage{times}  
\usepackage{graphicx} 

\date{}

\begin{document}


\title{Optimal Oil Production under Mean Reverting L\'{e}vy Models with Regime Switching}
\author{Moustapha Pemy \thanks{Department of Mathematics, Towson
University, Towson, MD 21252-0001, mpemy@towson.edu
 }
}
\maketitle


\begin{abstract} 
This paper is concerned with the problem of finding the optimal of extraction policies of an oil field  in light of various financial and economical restrictions and constraints. Taking into account the fact that the oil price in worldwide commodity markets fluctuates randomly following global and seasonal macro-economic parameters, we model the evolution of the oil price as a mean reverting regime switching jump diffusion process. We formulate this problem as finite-time horizon optimal control problem.  We solve the control problem using the method of viscosity solutions. Moreover, we construct and prove the convergence of  a numerical scheme for approximating the optimal reward function and the optimal extraction policy. A numerical example that illustrates these results is presented.  
 
\vskip 0.2 in

\end{abstract}





\section{Introduction}

Oil and natural gas have always been the main sources of revenue for a large number of developing countries as well as some industrialized nations around the world. Oil extraction policies vary from a country to another, in some countries the extraction is done by a state owned corporation in others it is done by foreign multinationals. The production and regulation of strategic natural resources such as oil, natural gas, uranium, gold, copper,...,etc  have always been among of the leading topics in geopolitical and macro-economical debates between politicians as well as  financial economists in academic circles.\\
The optimal extraction of natural resources was first studied done by Hotelling (1931), he derived  an optimal extraction policy under the assumption that the commodity price is constant. Many economists have extended the Hotelling model by taking into account the uncertainty in the supply, the demand,  as well as the ever-changing regulatory landscape 
of natural resource policies. Among many others, one can cite the work of Sweeney (1977),   Hanson (1980), Solow and Wan (1976), Pindyck (1978), (1980),  Lin and  Wagner (2007) for various extensions of the basic Hotelling model. Cherian {\it et al.} (1998) studied the optimal extraction of nonrenewable resources as a stochastic optimal control problem with two state variables, the commodity price and the size of the remaining reserve. They solved the control problem numerically by using Markov chain approximation methods developed by Kushner (1977) and Kushner and Dupuis (1992). Recently Aleksandrov {\it et al.} (2012) studied the optimal production of oil as an American-style real option and used Monte-Carlo methods to approximate the optimal production rate when the oil price follows a mean-reverting process.\\
It is also self evident that the price of oil in commodity exchange markets  fluctuates following divers macro economical and global  geopolitical forces. It is therefore crucial to take into account the random dynamic of the commodity value when solving the optimal extraction problem in order to maintain  the validity of the result obtained. In this paper, we use the mean reverting regime switching L\'{e}vy processes to model the oil price. Mean reverting processes were first used to describe the evolution of commodity prices  by Gibson an Schwartz (1990) and   Schwartz (1997), these processes capture perfectly the mean reversion feature of commodity prices around an equilibrium price when the market is stable. Other authors like Schwartz and Smith (2000) and Aleksandrov {\it et al.} (2012) used the mean reverting processes in their extraction models. Oil prices also display a great deal of seasonality, jumps and spikes due to various supply disruptions and political turmoils in oil-rich countries,  we use regime switching jump diffusions to capture all those effects. Thus our pricing model closely captures the instability of oil markets. Regime switching models have been extensively used in the financial economics literature since its introduction by Hamilton (1989).  Many 
authors have studied the control of systems that involve regime switching
 using a hidden Markov chain, one can cite Zhang and Yin (1998), (2005), Pemy and Zhang (2006), Pemy (2011), (2014) among others. \\ 
In this paper we treat the problem of finding optimal extraction strategies as an optimal control problem of a mean reverting Markov switching L\'{e}vy processes in a finite time horizon. The main contribution of this paper is two-fold, first we prove that the value function is the  unique viscosity solution of the associated 
Hamilton-Jacobi-Bellman equation. Then, we build a finite difference approximation scheme  and prove its convergence to the unique viscosity solution of HJB equation. This enable us to derive both the optimal  reward function and the optimal extraction policy in this broad setting.\\ 
The paper is organized as follows. In the next section, we
formulate the problem under consideration. In section 3 we derive the properties of the value function and characterize it as unique viscosity solution of the Hamilton Jacobi Bellman equation. And in section 4, we construct  a finite difference approximation scheme and prove its convergence to the value function. Finally, in section 5, we give a numerical example.
\section{Problem formulation}
Consider a multinational oil company with an extraction lease that expires in $T$ years,   $0<T<\infty$. We assume that  the market value of a barrel of oil at time $t$ is $S_t=e^{X_t}$. Given that oil prices are very sensitive global macro-economical and geopolitical shocks, we model $X_t$ as a  mean reverting regime switching L\'{e}vy process with two states. Let $\al(t)\in{\cal M}=\{1,2\}$ be a finite state Markov chain that captures the state of the oil market: $\al(t)= 1$ indicates the bull market at time $t$ and $\al(t) = 2$ represents a bear market at  time $t$.  The generator of this Markov chain is 
\bea
Q=\left(\begin{array}{ll} -\lambda_1&\lambda_1\\
\lambda_2&-\lambda_2 \end{array} \right),\qquad \lambda_1>0,\lambda_2>0.
\eea
Let  $(\eta_t)_t$ be a L\'{e}vy process and let $N$ be the Poisson random measure of $(\eta_t)_t$, $  N(t,U)=\sum_{0<s\leq t}{\bf 1}_U(\eta_s-\eta_{s^-})$  for any Borel set $U\subset \rr$. Moreover, let $\nu$ be the L\'{e}vy measure of $(\eta_t)_t$ we have $\nu(U)=E[N(1,U)]$ for any Borel set $U\subset \rr$. The differential form of $N$ is denoted by $N(dt,dz)$, we define the  differential $\bar{N}(dt,dz)$ as follows
\bea 
\disp
\bar{N}(dt,dz)=\left\{\begin{array}{ll} N(dt,dz)-\nu(dz)dt \qquad &\hbox{if  } |z|<1\\
   N(dt,dz)&\hbox{if  } |z|\geq 1.   \end{array} \right.
   \eea
We assume that  the L\'{e}vy measure $\nu$ has finite intensity,
\beq{finiteIntensity} 
\disp \Gamma = \int_{\rr} \nu(dz)<\infty.
\eeq
 In other terms, the total sum of jumps and spikes of the oil price during the lifetime of the contract is finite.
Let $ K<\infty$ be the total size of the oil field at the beginning of the lease, and let $Y(t)$ be the size of the remaining reserve of the oil field by time $t$, obviously $Y(t)\in [0,K]$. The state variables of our control problem are $X(t)\in\rr$ and $Y(t)\in [0,K]$, and the state space is $[0, \infty)\times[0,K]$. 
 We  assume that  the processes $ X(t), Y(t)$ follow the   dynamics
\beq{model}
 \left \{ \begin{array}{ll}
 \disp\mathrm{d}X(t)&\disp =\kappa \big(\mu(\al(t))- X(t)\big)\mathrm{d}t+\sigma(\alpha(t))\mathrm{d}W(t) + \int_{\rr}\gamma(\al(t))z\bar{N}(dt,dz),\\
 dY(t)&= -u(t) dt, \\
\disp X(s)&=x\geq 0, \,\,\,\, Y(s) = y \geq 0, \qquad 0\leq s\leq t\leq T, \end{array} \right.
\eeq   
 where $ u(t)\in U=[0,\bar{u}]$ is the extraction rate chosen by the company. In fact the process $ u(t)$  is control variable. The process  $ W(t)$ is the Wiener process  defined on a probability space  $ (\Omega,{\cal F },P)$. Moreover, we assume that $W(t)$, $\al(t)$ and $\eta_t$  are independent. The parameter $\mu(\cdot)$ represents the equilibrium price of oil. For each state $i\in \{1,2\}$ of the oil market we assume that the corresponding equilibrium price $\mu(i)$ is known. Similarly $\sigma(\cdot)$ represents the volatility and $\gamma(\cdot)$ represents the intensity of the jump diffusion. For each state $i\in\{1,2\}$ of the oil market  we assume that $\sigma(i),\, \gamma(i)$ are known constants. As a matter of fact $\gamma(i)$ captures the frequencies and sizes  of jumps and spikes of the oil price. It is well known that  the spot prices of energy commodities that are very expensive to store usually have frequent jumps and spikes within  short periods of time. In order to capture those effects jump diffusions are the ideal candidates.  \
\begin{defn}
The extraction rate  $u(\cdot)$  taking values on intervals $[0,\bar{u}]$  is called  an admissible control with respect to the initial data $(s,x,y, i)\in[0,T]\times[0,\infty)\times[0,K]\times{\cal M}$ if:
 \begin{itemize}
\item Equation (\ref{model}) has a unique  solution  with $X(s) = x$,  $Y(s) =y$, $\al(s)= i $, and $X(t)\in[0,\infty), \,\, Y(t)\in[0,K]$ for all $t\in [0,T]$.
\item The process $u(\cdot)$  is $\{{\cal F}_t\}_{t\geq 0}$-adapted where $ {\cal F}_t = \sigma\{\alpha(s),W(s), \eta_s; s\leq t\}$.
\end{itemize}
We use ${\cal U} ={\cal U}(s,x,y,i)$ to denote the set of admissible controls taking values in $U=[0,\bar{u}]$ such that $X(s)=x$, $Y(s) =y$, $\al(s)=i$.  
\end{defn}
 The admissibility condition implies that we should only consider extraction rates that depend on the information available up to time $t$, and within a reasonable range sets forth  at the beginning of the lease and that will guarantee that the state variables  stay in state space during the lifetime of the lease. 
We assume that  the extraction cost function per unit of time $t$ is the function $C(t,Y(t),u(t))$ that depends on the extraction rate $u(t)$ and the size of the remaining reserve $Y(t)$. Moreover,  $C(t,y, u)$ should be measurable  and nondecreasing with respect to size of the remaining reserve $y$. This enables us to capture  the fact that is more expensive to extract as the size of the oil field decreases.  A typical example of an extraction cost function is 
\bea
\disp C(t,y,u) = a+ m u(b Y+c),
\eea
where $a>0$ can be seen as the initial cost  of setting up the mine and $m$,  $b$, and $c$ are constants such that $m>0$ and $b,c\geq 0$.
The total profit rate for operating the mine is
\bea
 L(t,X(t),Y(t),u(t))=e^{X(t)} u(t) -C(t,Y(t), u(t))\quad t\in[0,T).
\eea
We assume at the end of the lease there are no revenues from extraction but the oil field still has some value and we roughly estimate that to be equal to overall market value of the remaining oil under the ground. We model that terminal value as follows
 \[
\Psi(T,X(T),Y(T)) =(K-Y(T))(e^{X(T)} -m) .
\]
Given a discount rate $r> 0$, the payoff functional  is 
\beq{funcPlayer1}
&& J(s,x,y,\iota;u)\nonumber\\
&&= E\Bigg[\int_s^Te^{-r(t-s)} L(t,X(t),Y(t),u(t) \rm{d}t  \nonumber  \\
&& \hspace{0.5in}+ e^{-r(T-s)}\Psi(T,X(T),Y(T))\bigg{|} X(s) = x, Y(s) = y, \al(s) = \iota\Bigg].
\eeq
The oil  company will try to maximize its payoff by adjusting the extraction rate $ u(\cdot)$ accordingly, the optimal reward function also known as the value function of the control problem is
\beq{val1}
 V(s,x,y,\iota) = \sup_{u\in{\cal  U}}J(s,x,y,\iota; u)=J(s,x,y,\iota; u^*). 
\eeq
  We define the Hamiltonian of the control problem as follows
\beq{Hamilton1}
&& H(s,x,y,\iota,V,\frac{\partial V}{\partial x},\frac{\partial V}{\partial y}, \frac{\partial^2 V}{\partial x^2}\bigg)     \nonumber\\
&&=rV - \sup_{u\in U} \Bigg{(}\frac{1}{2}\sigma^2(\iota)\frac{\partial^2V}{\partial x^2}  
+  \kappa \big(\mu(\iota)- x\big)\frac{\partial V}{\partial x}
 -  u\frac{\partial V}{\partial y}   + \nonumber  \int_{\rr}\bigg(V(s,x+\gamma(\iota)xz,y,\iota) -V(s,x,y,\iota) \nonumber \\
&&  -{\bf 1}_{\{|z|<1\}}(z)\frac{\partial V}{\partial x} \gamma(\iota)xz\bigg)\nu(dz)  + L(s,x,y,u,\iota)+  QV(s,x,y,\cdot)(\iota)\Bigg{)},
\eeq	 
with $\disp  QV(s,x,y,\cdot)(\iota) = \sum_{j\neq \iota}q_{\iota j}(V(s,x,y,j) -V(s,x,y,\iota))$. 
 In order  to find the optimal extraction strategy $ u^*$  we first have to derive the value function $ V$ of the control problem. 
  Formally the value function $ V$ should satisfy  the following Hamilton Jacobi Bellman equation.
\beq{isaac1}
 \left \{ \begin{array}{ll} \frac{\partial V}{\partial s}  =   H(s,x,y,\iota,V,\frac{\partial V}{\partial x},\frac{\partial V}{\partial y}, \frac{\partial^2 V}{\partial x^2}\bigg),& \quad (s,x,y,\iota)\in  [0,T)\times \rr^+\times [0,K] \times {\cal M} \\
V(T,x,y,\iota) = \Psi(T,x,y), &\quad (x,y,\iota)\in  \rr^+\times [0,K] \times {\cal M}  
\end{array}\right.
\eeq
In the next section we analyze the properties of the value function and fully characterize the value function as the unique viscosity solution of the Bellman equation (\ref{isaac1}).

\section{Properties of the Value Function}
The Bellman equation (\ref{isaac1}) is a system of coupled fully nonlinear integro-differential partial differential equations which may not have smooth solutions in general.  In order to solve  (\ref{isaac1}) we will use a weaker form of solutions namely the notion of viscosity solutions  introduced by Crandall and Lions (1983).   Let us first recall the definition of viscosity solution.  
\begin{defn}
The function $ W$ defined on $ {\cal D}:= [0,T]\times \rr^+\times [0,K] \times {\cal M} $  is a viscosity subsolution (resp.  supersolution ) of 
\beq{visco}
  \frac{\partial W}{\partial s}= H(s,x,y,\iota,W, \frac{\partial W}{\partial x}, \frac{\partial W}{\partial y},\frac{\partial^2 W}{\partial x^2}),
\eeq
 if $ W$ is lower semi-continuous (resp.  upper semi-continuous), and for any $ \iota \in {\cal M }$,  for any test function $ \phi\in {\cal C}^{1,2,1}( [0,T]\times \rr^+\times [0,K]$  such that $ W-\phi$ has a local maximum (resp. minimum) at $ ( s_0,x_0,y_0,\iota) \in {\cal D}$  
\beq{viscodef1}
   \frac{\partial \phi}{\partial s}(s_0,x_0,y_0,\iota)\leq H(s_0,x_0,y_0,\iota, W,\frac{\partial \phi}{\partial x}, \frac{\partial \phi}{\partial y},\frac{\partial^2 \phi}{\partial x^2}),
\eeq
\beq{viscodef2}
\bigg(\rm{resp.} \ad   \frac{\partial \phi}{\partial s}(s_0,x_0,y_0,\iota)\geq H(s_0,x_0,y_0,\iota, W,\frac{\partial \phi}{\partial x}, \frac{\partial \phi}{\partial y},\frac{\partial^2 \phi}{\partial x^2}) \bigg).
\eeq
$W$ is a viscosity solution of (\ref{visco}) if $W$ is both a viscosity subsolution and supersolution.
\end{defn}
\begin{lem}\label{conti}
For each $\iota\in{\cal M}$,
the value function $V(s,x,y,\iota)$ is Lipschitz continuous with respect to $s,x$ and $y$ and  has  at most a linear growth rate, i.e.,
 there exists a constant $C$ such that
$\mid V(s,x,y,\iota) \mid \leq C(1+|x|+|y|)$.
\end{lem} 
The continuity of the value function $ V$  naturally comes from the application of the It\^{o}-L\'{e}vy isometry, the Lipschitz continuity of the parameters of the model and the  Gronwall inequality. For more details, one can refer to Pemy (2014) for the proof of a similar result in the case of the optimal stopping of regime switching L\'{e}vy processes.

\begin{rem}
The dynamical Programming Principle  implies that
\beq{dyna1}
&& \disp\quad   V(s,x,y,\iota) \nonumber\\
  &=& \sup_{u\in{\cal  U}} E\bigg[\int_s^Te^{-r(t-s)}L(t,X(t),Y(t),u(t),\al(t)\rm{d}t  \\
  && +e^{-r(T-s)}V(T,X(T),Y(T),\al(T))\Bigg{|} X(s) = x, Y(s) = y, \al(s) = \iota\bigg]\nonumber 
\eeq
\end{rem}
Using the Dynamical Programing Principle and the continuity of the value function $V$   we can now characterize the value function $V$  as unique viscosity solution of the Hamilton Jacoby Bellman equation (\ref{isaac1}).
\begin{thm}\label{viscosity}
The value function $ V$ is the unique viscosity solution of the Bellman equation (\ref{isaac1})
\beq{visco1}
  \frac{\partial W}{\partial s}= H(s,x,y,\iota,W, \frac{\partial W}{\partial x}, \frac{\partial W}{\partial y},\frac{\partial^2 W}{\partial x^2}).
\eeq
\end{thm}
This result is a standard result in control  theory. For more about the viscosity solution characterization of the value function  one can refer to Fleming and Soner (1993),  $\O$ksendal and Sulem (2004), Pemy (2014)
 among others. For more on the theory and application of viscosity solutions 
one can refer to Crandall, Ishii and Lions (1992), Yong and Zhou (1999).   
The next result gives the road map we will use to find the optimal extraction strategy if we already have the value function.
\begin{thm}\label{NashEquilibrium}
Assume  that 
the nonlinear Hamilton Jacobi Bellman equation (\ref{isaac1})  has a solution $ V(s,x,y,\iota)$, let $u^*\in {\cal U}$ such that
 \beq{optControl1} 
  \disp  u^*(s) = \arg\max\bigg(   -  u\frac{\partial V}{\partial y}    +  L(s,x,y,u,\iota)  \bigg), \quad s\in[0,T].
\eeq
Then the $ u^*$ is the optimal strategy and 
$  J(s,x,y,\iota;u^*)  = V(s,x,y,\iota)$.
\end{thm}
This result is just the standard Verification Theorem in control theory, one can refer to Fleming and Rishel (1975)  and Fleming and Soner (1993) for more details. 
\section{Numerical Approximation}
In this section, we construct a finite difference scheme and
show that it converges to the unique viscosity solution of the Bellman equation (\ref{isaac1}). 
Let $k, h, l\in(0,1) $ be the step size with respect to $s$, $x$ and $y$ respectively,  we consider the finite difference
operators  $\Delta_s$, $\Delta_x$, $\Delta_{xx}$ and $\Delta_y$ defined by
\bea \ad \Delta_s V(s,x,y,i)=\frac{V(s+k,x,y,i)-V(s,x,y,i)}{k}, \quad \Delta_x V(s,x,y,i)=\frac{V(s,x+h,y,i)-V(s,x,y,i)}{h} \\
\ad \Delta_y V(s,x,y,i) =\frac{V(s,x,y+l,i)-V(s,x,y,i)}{l},\\ 
\ad \Delta_{xx} V(s,x,y,i)=\frac{V(s,x+h,y,i)+ V(s,x-h,y,i)-2V(s,x,y,i)}{h^2}.
 \eea
Let ${ \rm I}f$ denote the integral part of the Hamiltonian $H$. We will approximate ${\rm I}f$ using the Simpson quadrature. In fact we have 
\bea
\ad{\rm I}f(s,x,y,i)\\
\ad = \int_{\rr}\bigg(f(s,x+\gamma(i)zx,y,i)-f(s,x,y,i) -{\bf 1}_{\{|z|<1\}}(z)\frac{\partial f(s,x,y,i)}{\partial x}\cdot \gamma(i)zx\bigg)\nu(dz).
\eea
Using the fact the L\'{e}vy measure is finite  $\disp \Gamma = \int_{\rr} \nu(dz)<\infty$, we have 
\bea
\ad{\rm I}f(s,x,y,i)= \int_{\rr}f(s,x+\gamma(i)zx,y,i)\nu(dz)-\frac{\partial f(s,x,y,i)}{\partial x}\int_{-1}^1 \gamma(i)xz\nu(dz)  -f(s,x,y,i)\Gamma. 
\eea
We use the Simpson's quadrature to approximate the integral part of the Hamiltonian.
Let $\xi\in(0,1)$ be the step size of the Simpson's quadrature, the corresponding approximation of the integral part is
\bea
\ad{\rm I}_{\xi}f(s,x,y,i) = \sum_{j=0}^{N_\xi}c_jf(s,x+\gamma(i)xz_j,y,i)-\frac{\partial f(s,x,y,i)}{\partial x}\sum_{j=0}^{M_\xi} d_j\gamma(i)xz_j-f(s,x,y,i)\Gamma, 
\eea
where  $(c_j)_{0\leq j\leq N_\xi}$ and $(d_j)_{0\leq j\leq M_\xi}$ are the corresponding sequences of the coefficients of the Simpson's quadrature. In fact $\disp  \lim_{N_\xi \rightarrow \infty}  \sum_{j=0}^{N_\xi}c_j =\Gamma$ and $\disp \lim_{M_\xi \rightarrow \infty}\sum_{j=0}^{M_\xi}d_j=\int_{-1}^1\nu(dz) $.
The corresponding discrete version of the Hamiltonian $H$  is given by
\beq{disc-Ham1}
&&  H_{h,k,l}V(s,x,y,i) \nonumber \\
 &=& rV(s,x,y,i)-\sup_{u\in U} \Bigg(\frac{1}{2}\sigma^2(i)\Delta_{xx}V(s,x,y,i) + \kappa( \mu(i)-x)\Delta_x V(s,x,y,i) \\
&&- u\Delta_yV(s,x,y,i)+{\rm I}_{\xi}V(s,x,y,i) +L(s,x,y,u,i)
  + QV(s,x,y,\cdot)(i)\Bigg). \nonumber
\eeq
Therefore the discrete version of (\ref{isaac1}) is 
\beq{discrete}
 \left \{ \begin{array}{ll}\disp V(s,x,y,i)  = \frac{1}{r}\Delta_sV  +\frac{1}{r}H_{j,k,l}V(s,x,y,i).& \\
V(T,x,y,\iota) = \Psi(T,x,y).
\end{array}\right.
\eeq
First we prove the  existence of a solution for the discretized equation (\ref{discrete}) on bounded subsets of the domain of study ${\cal D}$ where  ${\cal  D}:=[0,T]\times\rr^+ \times [0,K] \times {\cal M}$. We define ${\cal D}_R=\{(s,x,y,i)\in {\cal D}; |x|<R\}$, for some $R>0$. We will restrict our study on the set ${\cal D}_R $ for some $R$ large enough. As a matter of fact,  we are just assuming that the oil price will not go beyond a reasonable large threshold.   We will approximate our solution on that bounded domain. We have the following crucial Lemma.  
\begin{lem}\label{DiscreteLemma}
Let  $\xi>0 $ be small enough, for each $h,k,l\in (0,1)$, there exists a  unique bounded function $V_{l,h,k}$  defined on ${\cal D}_R$  that solves equation (\ref{discrete}).
\end{lem}
\begin{rem}
\begin{enumerate}
\item Define  $S\to (0,1)^4\times[0,T]\times\rr^+\times[0,K]\times {\cal M}\times \rr\times B([0,T]\times\rr^+\times[0,K]\times {\cal M} )$ as follows;
\beq{schema}
&&S(\xi, h,k,l, s,x,y,i,w,W)\nonumber \\
&=& w -\sup_{u\in U}\Bigg( a(x,i)W(s,x+h,y,i) \nonumber \\
&&+ b(x,i)W(s,x-h,y,i) -c(x,i;u)w  - \frac{yu}{rl}W(s,x,y+l,i) \nonumber\\
&&
+ \sum_{j=0}^{N_\xi}\frac{c_j}{r}W(s,x+\gamma(i)xz_j,y,i)+ \sum_{j\ne i} \frac{q_{ij}}{r}(u)W(s,x,y,j)\nonumber\\
&&- W(s,x+h,y,i)\frac{ \sum_{j=0}^{M_\xi} d_j\gamma(i)xz_j}{rh}+ L(s,x,y,u,i)\Bigg),
\eeq
where coefficients  $c(x,i;u)$, $a(x,i)$  and, $b(x,i)$ are defined in (\ref{coefficient1}), (\ref{coefficient2}) and (\ref{coefficient3}) respectively.
Obviously $V_{h,k,l}$ solves the equation $S(\xi, h,k,l, s,x,y,i,V_{h,k,l}(s,x,y,i),V_{h,k,l}) =0$.
 It is clear that for $h$ small enough the coefficients $a(x,i)>0$, $b(x,i)>0$ therefore the scheme $S$ is monotone with respect to argument $W$ i.e., for  all $\xi,h,k,l \in(0,1),s\in [0,T], x\in \rr^+, y\in [0,K], i \in {\cal M}$ and $W_1, W_2 \in  B([0,T]\times\rr^+\times[0,K]\times {\cal M} )$, we have
\beq{monotone}
\quad S(\xi, h,k,l, s,x,y,i,w,W_2)\leq S(\xi, h,k,l, s,x,y,i,w,W_1) \quad \hbox{whenever} \quad W_1\leq W_2.
\eeq
\item 
It is clear from \lemref{DiscreteLemma} that the numerical scheme obtained from (\ref{discrete}) is stable  since the solution of the scheme is bounded independently of the step sizes $h,k, l\in(0,1)$ and obviously consistent because as the step sizes $h,k,l$ go to zero the finite difference operators converge to the actual partial differential operators. We have the following convergence theorem.

\end{enumerate}

\end{rem}

\begin{thm}\label{convergence}
For each $\xi>0$ small enough, let $V_{h,k,l}$ be the solution of the discrete scheme obtained in \lemref{DiscreteLemma}. Then as $\xi \downarrow 0$ and  $(h,k,l)\rightarrow 0$ the sequence $V_{h,k,l}$ converges locally uniformly on $ {\cal D}_R$ to the unique viscosity solution $V$ of (\ref{isaac1}).
\end{thm}
This result  is the standard method for approximating  viscosity solutions, for more one can refer to Barles and Souganidis (1991).
Below is the implementation algorithm.\\
{\it Fixed Point Algorithm}
\begin{enumerate}
\item Choose a tolerance $\epsilon>0$.
Choose an initial guess  of $V$ denoted by $V^{(0)}$
 \item For $j=0,...,\rm{MaxIteration}.$ \\
(a)  Find $u^*$ such that
 \beq{optControl3} 
  \disp  u^* = \arg\max_{u\in U}\bigg(   -  u\frac{\partial V^{(j)}}{\partial y}    +  L(s,x,y,u,\iota)  \bigg), 
\eeq
 
  (b) Solve the equation 
\bea
  V^{(j+1)}  = \frac{1}{r}\Delta_sV^{(j)}  +\frac{1}{r}H_{h,k,l}V^{(j)},
\eea
\item If $\|V^{(j+1)}-V^{(j)}\|<\epsilon$, then stop, else go to step 2 with $j\leftarrow j+1$. 
\end{enumerate}


\section{Applications}
The oil field has a known capacity of K=10 billion barrels and the lease has a T=10 years maturity.  The market has two trends: the up trend and the down trend.  When the market is up,  $ \al(t)=1 $, the oil equilibrium price  is $\mu(1)=55$ and when the market is down, $\al(t)=2$, the equilibrium price is  $\mu(2) = 35$. The mean reversion coefficient is $\kappa = 0.01$,
the volatility   is $\sigma(1) = 0.2$ when the market is bullish and  $ \sigma(2) =0.3$ when the market is bearish.
 And the jump intensity is $ \lambda(1) =0.01$ when the market is up and $ \lambda(2)=0.15$ when the market is down. We assume that 
The profit rate function of the oil company per unit of time (hour) for each barrel of crude oil extracted is
\[
L(t,x,y,u) = (e^{x}u-(5+20 u)).
\] The  terminal value of the oil filed is
\[
\Phi(T,x,y) =(K-y)(e^{x}-20).
\]
Moreover, we assume that the extraction $ u(\cdot) \in [0,50000]$. Keep in mind that, because the payoff rates are linear functions of each control variable $ u(\cdot)$.
Once the value function $ V$ is approximated  numerically,  using \thmref{NashEquilibrium} the optimal strategy $ u^*$  is obtained by looking at the sign of the derivative  of the quantity $  -u \frac{\partial V}{\partial y} + (xu-(5+20 u))$ with respect to $ u$. Let $G$  be that derivative, we have
\bea
G(s,x,y,i) = -\frac{\partial V}{\partial y} + (e^{x}-20).
\eea
 
 We see that  the optimal extraction strategy will only be attained at the endpoints of the intervals $ U=[0,50000] $, we have.
\bea 
 u^*(s) = \left \{ \begin{array}{ll} 0 & \hbox{ if }G(s,x,y,i)\leq 0\\
50000 &\hbox{ if }G(s,x,y,i)>0.\end{array}\right.
\eea
 In Figures 1 and 2, we have plots of the function $G(s,x,y,i)$. Note that the sign of this function will dictate our optimal  extraction policy. In all these plots, the region above the line represents the domain where it is always optimal to extract at full capacity and the region below the curve represents the domain where it is better not to extraction. This is a typical case of a bang-bang control that is easy to implement.  
\begin{figure}{fig1}
\includegraphics[height=20cm,width=20cm]{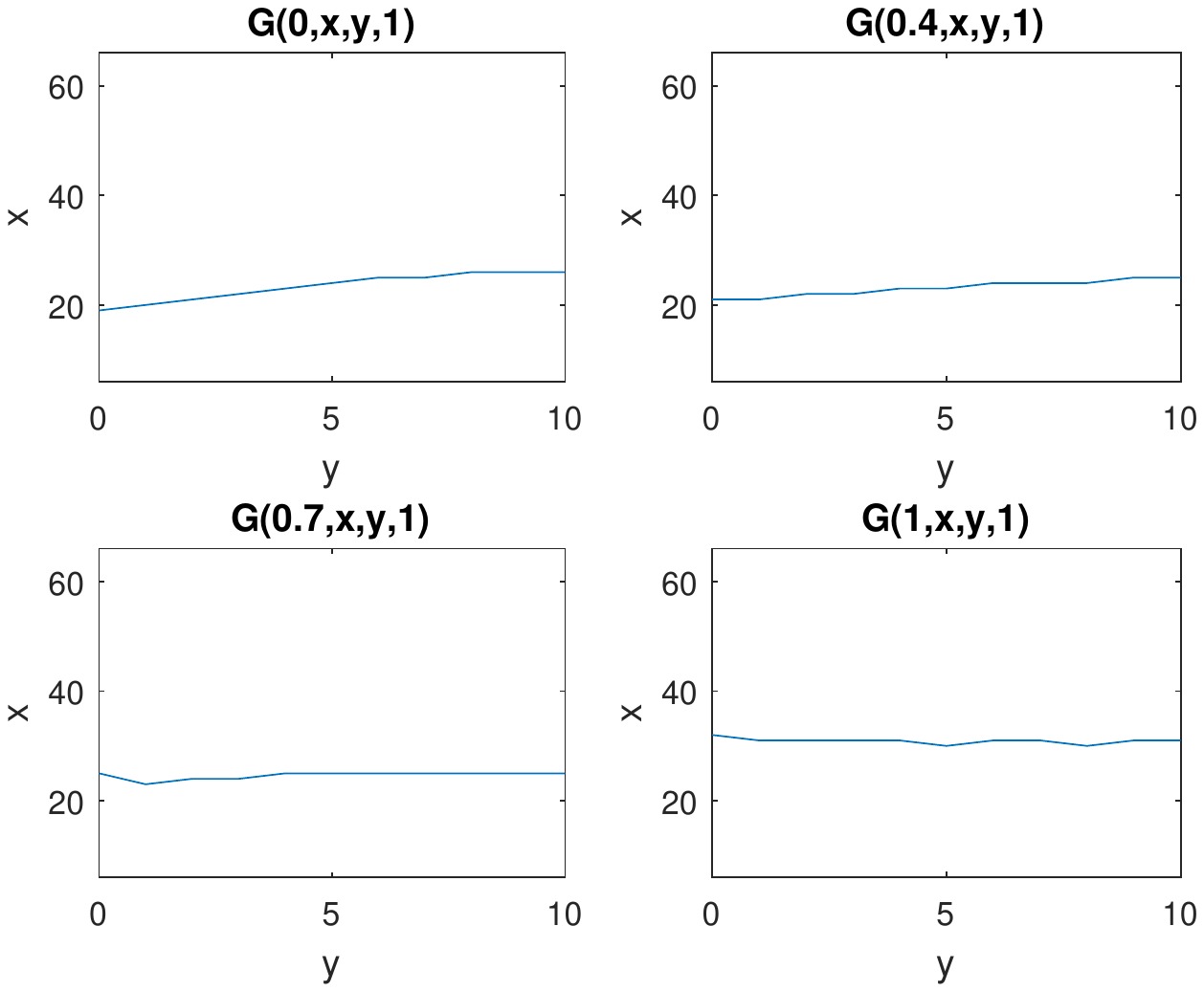}
\caption{Plots of $G(s,x,y,i)$ at various times $s\in \{0, 0.4, 0.7, 1\}$  when the market is up. }
\end{figure}

\begin{center}
\begin{figure}{fig2}
\includegraphics[height=20cm,width=20cm]{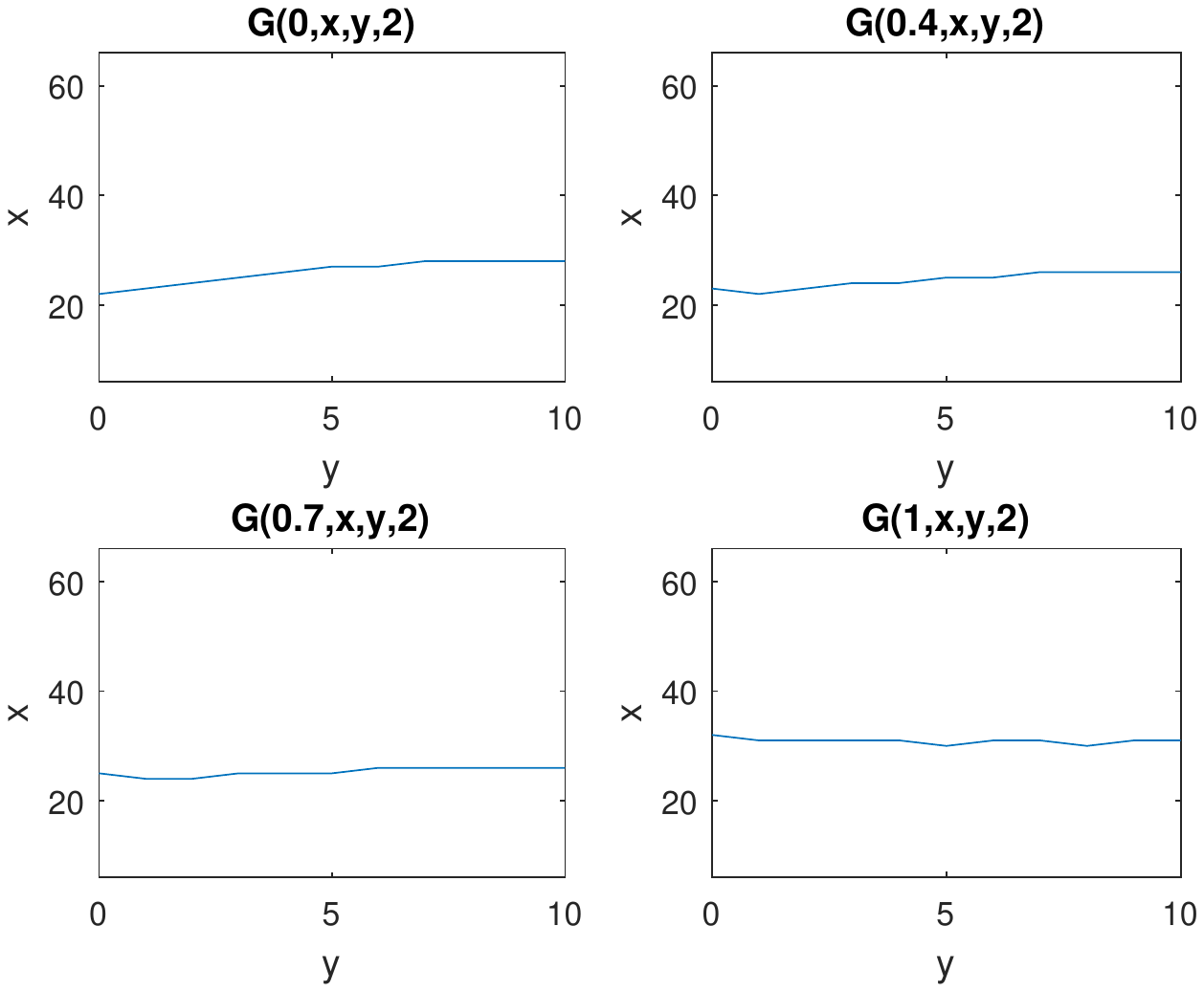}
\caption{Plots of $G(s,x,y,i)$ at various times $s\in\{0, 0.4, 0.7, 1\}$  when the market is down.}
\end{figure}
\end{center}
\appendix
\section{Appendix: Proofs of Results}
\subsection{Proof of \lemref{conti}}
In fact it can be shown that the value function are Lipschitz continuous with respect to $x$ and $y$. A detailed proof can be found in Pemy (2014)  in the particular case of an optimal stopping problem.    Below we show the linear growth property of the value function.
The linear growth inequality follows from the Lipschitz continuity of the value function with respect to $x$ and $y$. Thus  there exist $C, C'>0$ such that
\bea \disp
|V(s,x,y,i)|\leq C'|x| +|V(s,0,y,i)|,
\eea
and
\bea
\disp
|V(s,0,y,i)|\leq C|y|+|V(s,0,0,i)|.
\eea
Combining the last two inequalities gives,
\bea
\disp
|V(s,x,y,i)|\leq \max(C',C)(|x|+|y|+|V(s,0,0,i)|)\leq C_0(|x|+|y|+1),
\eea
for some $C_0>\max(C',C)$.
This completes the proof. 
\endproof

\subsection{Proof of  \thmref{viscosity}}
Let $\iota \in {\cal M}$, and $\psi \in \mathnormal{C}^{1,2,1}([0,T]\times\rr^+\times[0,K])$
such that $V(s,x,y,\iota)-\psi(s,x,y)$ has local minimum at $(s_0,x_0,y_0)$ in a neighborhood $N(s_0,x_0,y_0)$. Without loss of generality we assume that $V(s_0,x_0,y_0,\alpha_0)-\psi(s_0,x_0,y_0).$
We set $\al(s_0) = \al_0$ and  define a function 
\beq{varphi2}
\varphi(s,x,y,\iota)=\left \{ \begin{array}{ll}
                    \psi(s,x,y) &  \textrm{if} \qquad i=\alpha_0,\\
V(s,x,y,\iota),  & \textrm{if} \qquad i \not =\alpha_0. \end{array} \right. 
\eeq
Let $\gamma\geq s_0 $ be the first jump time of $\alpha(\cdot)$ from the
initial state $\alpha(s_0)=\al_0$, and let $\eta \in [s_0,\gamma]$ be such that 
$(t,X(t),Y(t))$ starts at $(s_0,x_0,y_0)$ and stays in $N(s_0,x_0,y_0) $ for $s_0\leq t \leq \eta $.
Moreover, $ \alpha(t)= \alpha_0 $, for $ s_0\leq t\leq \eta$.
Let $u(\cdot)$ be an admissible control such that $u(t)= u$ for $t\in[0, \eta]$. 
From then Dynamical Programming Principle (\ref{dyna1}) we derive
\beq{inesup}
V(s_0,x_0,y_0,\al_0)\geq E\bigg[\int_{s_0}^{\eta}e^{-r(t-s_0)} L(t,X(t),Y(t),u(t),\al(t))dt  \nonumber \\
+ e^{-r(\eta-s_0)}V(\eta,X(\eta),Y(\eta),\al(\eta))\bigg{]}.
\eeq
Using Dynkin's formula we have,
\beq{dynkin2}
& & E^{s_0,x_0,y_0,\al_0}[ e^{-r(\eta-s_0)} \varphi(\eta,X(\eta),Y(\eta),\alpha_0)]- \varphi(s_0,x_0,y_0,\alpha_0) \\
&=&E^{s_0,x_0,y_0,\al_0} \int_{s_0}^\eta e^{-r(t-s_0)}[-r\varphi(t,X(t),Y(t),\al_0)+({\cal L}^{u}\varphi)(t,X(t),Y(t),\al_0)]dt.\nonumber
\eeq
where ${\cal L}^{u}$ is the generator of the Markov processes $(X_t, Y_t)$.
 Note that ${\cal L}^{u}(\varphi)(s,x,y,\iota)$ can be written as ${\cal L}^{u}(\varphi)(s,x,y,\iota)={\cal A}^{\iota,u}(\psi)(s,x,y)+ Q\varphi(s,x,y,\cdot)(\iota)$ with 
\bea
{\cal A}^{\iota,u}(\psi)(s,x,y) \aad= \frac{\partial \psi}{\partial s}+ \frac{1}{2}\sigma^2(\iota)\frac{\partial^2\psi}{\partial x^2}  
+ \kappa(\mu(i)-x)\frac{\partial \psi}{\partial x}
 -  u\frac{\partial \psi}{\partial y}   + \nonumber  \int_{\rr}\bigg(\psi(s,x+\gamma(\iota)xz,y) \\
\aad-\psi(s,x,y)   -{\bf 1}_{\{|z|<1\}}(z)\frac{\partial \psi(s,x,y)}{\partial x}\cdot \gamma(\iota)xz\bigg)\nu(dz). 
\eea
Given that $(s_0,x_0,y_0)$ is the minimum of $V(t,x,y,\alpha_0)-\psi(t,x,y)$ in $N(s_0,x_0,y_0)$.
For $ s_0\leq t \leq \eta$, we have
\beq{vbig2}
 V(t,X(t),Y(t),\alpha_0) \geq \psi(t,X(t),Y(t)) +V(s_0,x_0,y_0,\alpha_0) \nonumber \\
-\psi(s_0,x_0,y_0)=\varphi(t,X(t),Y(t),\alpha_0).
\eeq 
Using equation (\ref{varphi2}) and (\ref{vbig2}), we have  
\beq{dynkin3}
&  & E^{s_0,x_0,y_0,\al_0}[ e^{-r(\eta-s_0)}V(\eta,X(\eta),Y(\eta),\alpha_0)]- V(s_0,x_0,y_0,\alpha_0)\nonumber \\
&\geq& E^{s_0,x_0,y_0,\al_0} \int_{s_0}^\eta e^{-r(t-s_0)} \bigg[ {\cal A}^{\al_{s_0},u}(\psi)(t,X(t),Y(t))+Q\varphi(t,X(t),Y(t),\cdot)(\al_0)\nonumber \\
& &-rV(t,X(t),Y(t),\al_0)\bigg]dt. 
\eeq
Moreover, we have 
\beq{Qarray2}
Q\varphi(t,X(t),Y(t),\cdot)(\alpha_0) &=& \sum_{\beta \not = \alpha_0}q_{\alpha_0 \beta}\Bigg{(}\varphi(t,X(t),Y(t),\beta)-\varphi(t,X(t),Y(t),\alpha_0)\bigg{)}\nonumber \\
&\geq& \sum_{\beta \not = \alpha_0}q_{\al_0 \beta}\Bigg{(}V_1(t,X(t),Y(t),\beta)-V(t,X(t),Y(t),\alpha_0)\bigg{)}\nonumber\\
&\geq&QV_1(t,X(t),Y(t),\cdot)(\alpha_0).
\eeq
It follows from (\ref{inesup}), (\ref{dynkin3}) and (\ref{Qarray2}) that
\bea
\disp
E^{s_0,x_0,y_0,\al_0} \int_{s_0}^\eta e^{-r(t-s_0)} \bigg[ {\cal A}^{\al_{0},u,}(\psi)(t,X(t),Y(t))+QV(t,X(t),Y(t),\cdot)(\al_0)\\-rV(t,X(t),Y(t),\al_0)+L(t,X(t),Y(t),u(t),\al(t))\bigg]dt \leq 0.
\eea
Dividing by $\eta -s_0>0$ and then sending  $\eta \rightarrow s_0$ leads to
\beq{super-ineq}
-rV(s_0,x_0,y_0,\al_0)+{\cal A}^{\al_{0},u}(\psi)(s_0,x_0,y_0)\nonumber \\ 
+QV(s_0,x_0,y_0,\cdot)(\al_0)+L(s_0,x_0,y_0,u,\al_0)\leq 0.
\eeq
Since this inequality is true for any arbitrary control $u(t)\equiv u\in [0,\bar{u}]$, then taking the supremum over all values $u\in U= [0,\bar{u}]$ we have
\beq{Supersol1}
\disp rV(s_0,x_0,y_0,\al_0) -\sup_{u\in U}\Bigg({\cal A}^{\al_0,u}(\psi)(s_0,x_0,y_0)\nonumber \\
+QV(s_0,x_0,y_0,\cdot)(\al_0)+L(s_0,x_0,y_0,u,\al_0)\Bigg)\geq 0.
\eeq
The inequalities (\ref{Supersol1})  obviously proves  that the value function $V$ is a viscosity  supersolution  as defined in (\ref{viscodef2}).\\
Now, let us prove the subsolution inequality (\ref{viscodef1}). We want to show that for each $\iota \in {\cal M}$, 
\beq{rrmax}
rV(s_0,x_0,y_0,\iota) -\sup_{u\in{ U}}\Bigg({\cal A}^{\al_0,u}(\psi)(s_0,x_0,y_0)\nonumber \\
+QV(s_0,x_0,y_0,\cdot)(\iota)+L(s_0,x_0,y_0,u,\iota)\Bigg) \leq 0,
\eeq 
where $(s_0,x_0,y_0)$ is a local maximum of $V(s,x,y,\iota)-\psi(s,x,y)$.
Let us assume otherwise that the inequality (\ref{rrmax}) does not hold. In other terms, we assume that we can find a state $\al_0 \in {\cal M}$, values $(s_0,x_0,y_0)$ and a function $\phi\in {\cal C}^{1,2,1}([0,T]\times \rr^+ \times [0,K])$ such that $V(t,x,y,\alpha_0)-\phi(t,x,y)$ has a local maximum at $ (s_0,x_0,y_0)\in [s,T]\times \rr\times\rr^+,$  and we have
\beq{absurd}
rV(s_0,x_0,y_0,\al_0) -\sup_{u\in{ U}_1}\Bigg({\cal A}^{\al_0,u}(\psi)(s_0,x_0,y_0)\nonumber \\
+QV(s_0,x_0,y_0,\cdot)(\al_0)+L(s_0,x_0,y_0,u,\al_0)\Bigg)\geq \delta.
\eeq
for some constant $\delta>0$.\\
Let us assume without loss of generality that $V(s_0,x_0,y_0,\alpha_0)-\phi(s_0,x_0,y_0)= 0$.
We define
\beq{varpsi10}
\varphi(s,x,y,i)=\left \{ \begin{array}{ll}
                    \phi(s,x,y), &  \textrm{if} \qquad i=\alpha_0,\\
V(s,x,y,i),  & \textrm{if} \qquad i \not =\alpha_0 .\end{array} \right. 
\eeq
Let $\gamma $ be the first jump time of $\alpha(\cdot)$ from the state $\alpha_0$, and let $\eta_0 \in [s_0,\gamma]$ be such that $(t,X(t),Y(t))$ starts at $(s_0,x_0,y_0)$ and stays in $N(s_0,x_0,y_0) $ for $s_0\leq t \leq \eta_0 $.
Since $\theta_0 \leq \gamma $ we have $ \alpha(t)= \alpha_0 $, for $ s_0\leq t\leq \eta_0$.
Moreover, since $V(s_0,x_0,y_0,\alpha_0)-\phi(s_0,x_0,y_0)= 0$ and attains its maximum at $(s_0,x_0,y_0)$ in $N(s_0,x_0,y_0)$ then 
\[
V(\eta,X(\eta),Y(\eta),\alpha(\eta))\leq \phi(\eta,X(\eta),Y(\eta))\quad \hbox{for  }\,\,\, s_0\leq \eta \leq \eta_0.
\]  
 Thus, we also have
 \beq{Phi}
V(\eta,X(\eta),Y(\eta),\alpha(\eta))\leq \varphi(\eta,X(\eta),Y(\eta),\alpha(\eta))\quad \hbox{for  }\,\,\, s_0\leq \eta \leq \eta_0.
\eeq
Using the Dynamical Programming Principle (\ref{dyna1}), it clear that for any admissible control $u(\cdot)$ and  time $\tau$ such that $s_0<\tau \leq \eta_0$, we have 
\bea \disp
 J_1(s_0,x_0,y_0,\al_{0};u)&\leq \disp E^{s_0,x_0,y_0,\al_0}\bigg[\int_{s_0}^{\tau}e^{-r(t-s_0)} L(t,X(t),Y(t),u(t),\al(t))dt \\
&  + e^{-r(\tau-s_0)}V_1(\tau,X(\tau),Y(\tau),\al(\tau))\bigg{]}\\
 & \disp \leq E^{s_0,x_0,y_0,\al_0}\bigg[\int_{s_0}^{\tau}e^{-r(t-s_0)} L_1(t,X(t),u(t),,\al(t))dt\\
& + e^{-r(\tau-s_0)}\varphi(\tau,X(\tau),Y(\tau),\al(\tau))\bigg{]}.
 \eea
 Note that
\beq{Qmat21}
Q\varphi(t,X(t),Y(t),\cdot)(\alpha_0) & =  & \sum_{\beta \not = \alpha_0}q_{ \alpha_0 \beta}(V_1(t,X(t),Y(t),\beta)-\phi(t,X(t),Y(t))) \nonumber \\
& \leq  & \sum_{\beta \not = \alpha_0}q_{ \alpha_0 \beta}(V_1(t,X(t),Y(t),\beta)- V_1(t,X(t),Y(t),\alpha_0))  \nonumber \\
& \leq & QV_1(t,X(t),Y(t),\cdot)(\alpha_0).
\eeq
 Using the inequality (\ref{absurd}) we have 
 \beq{new-ineq}
& &J(s_0,x_0,y_0,\al_0;u) \nonumber \\
&\leq& \disp E^{s_0,x_0,y_0,\al_0}\bigg(\int_{s_0}^{\tau}e^{-r(t-s_0)}\bigg{\{} -\delta +rV(t,X(t),Y(t),\al_0)\nonumber\\
&& -{\cal A}^{\al_0,u}(\phi)(t,X(t),Y(t))
-  QV(t,X(t),Y(t),\cdot)(\al_0) \bigg{\}}dt\disp\nonumber\\
&&+e^{-r(\tau-s_0)}\varphi(\tau,X(\tau),Y(\tau),\al_0)\bigg{)}.
\eeq
 The Dynkin's formula, (\ref{varpsi10}) and (\ref{Qmat21})   imply that  
\beq{ligne20}
& & E^{s_0,x_0,y_0,\al_0} e^{-r(\tau-s_0)}\varphi(\tau,X(\tau),Y(\tau),\alpha_0) \nonumber\\
&=& E^{s_0,x_0,y_0,\al_0} \int_{s_0}^\tau e^{-r(t-s_0)} \Bigg{[}{\cal A}^{\al_0,u,=}(\phi)(t,X(t),Y(t))+ Q\varphi(t,X(t),Y(t),\cdot)(\alpha_0) \nonumber \\
& &-r\varphi(t,X(t),Y(t),\alpha_0) \bigg{]}dt  +\varphi(s_0,x_0,y_0,\alpha_0) \nonumber\\
&\leq& E^{s_0,x_0,y_0,\al_0} \int_{s_0}^\tau e^{-r(t-s_0)} \Bigg{[}{\cal A}^{\al_0,u}(\phi)(t,X(t),Y(t))+ QV(t,X(t),Y(t),\cdot)(\alpha_0) \nonumber\\
 &&-rV(t,X(t),Y(t),\alpha_0)  \bigg{]}dt+V(s_0,x_0,y_0,\alpha_0).
\eeq
Combining (\ref{new-ineq}) and (\ref{ligne20}) we have
\beq{end}
J(s_0,x_0,y_0,\al_0;u)\leq  E^{s_0,x_0,y_0,\al_0}\bigg(-\int_{s_0}^{\tau}e^{-r(t-s_0)}\delta dt\bigg)+V(s_0,x_0,y_0,\alpha_0).
\eeq
It is easy to see that the quantity $\gamma=\disp E^{s_0,x_0,y_0,\al_0}\bigg(\int_{s_0}^{\tau}e^{-r(t-s_0)}\delta dt\bigg)>0$, thus taking the supremum over all admissible control $u(\cdot)\equiv u$ we obtain
\beq{finish}
V(s_0,x_0,y_0,\al_0)\leq -\gamma +V(s_0,x_0,y_0,\al_0),
\eeq
which is a contradiction. This proves that the inequality (\ref{rrmax}) is satisfied. Obviously we derive the subsolution inequality (\ref{viscodef1}). 
Therefore, $V$ is a viscosity solution of (\ref{isaac1}). The uniqueness of the viscosity solution follows from the standard Ishii method, for more one can refer to Pemy (2014) for a similar proof of the uniqueness of viscosity solution of optimal stopping problems for regime switching L\'{e}vy processes. 
\endproof

\subsection{Proof of \lemref{DiscreteLemma}}
We define the operator ${\cal F}_\xi$  on bounded functions on ${\cal D}_R$ as follows
\beq{DiscrOpe}
& & {\cal F}_\xi(V)( s,x,y,i; h,k,l)\nonumber \\
 &=&  \frac{1}{r}\Delta_sV  +\frac{1}{r}H_{u^*,}^{j,k,l}V(s,x,y,i)\nonumber \\
&=&  \frac{1}{rk}V(s+k,x,y,i)  +\sup_{u\in U}\Bigg( a(x,i)V(s,x+h,y,i) \nonumber \\
&&+ b(x,i)V(s,x-h,y,i) -c(x,i;u)V(s,x,y,i)  - \frac{u}{rl}V(s,x,y+l,i) \nonumber\\
&&
+ \sum_{j=0}^{N_\xi}\frac{c_j}{r} V(s,x+\gamma(i)z_jx,y,i)+\frac{1}{r} L(s,x,y,u,i)+ \sum_{j\ne i}\frac{q_{ij}}{r}V(s,x,y,j)\nonumber\\
&&- V(s,x+h,y,i)\frac{ \sum_{j=0}^{M_\xi} d_j\gamma(i)z_jx}{rh}\Bigg) \quad \hbox{if }(s,x,y,i)\in  {\cal D}_R \hbox{ and } s<T,\\
&&{\cal F}_\xi(V)(T,x,y,i;h,k,l) = \Psi(T,x,y).\nonumber
\eeq 
Where the coefficients $a(x,i), b(x,i)$ and $c(x,i;u)$ are defined as follows 
\beq{coefficient1}
\ad c(x,i;u) =\frac{1}{ rk} + \frac{\sigma^2(i)}{rh^2} +  \frac{\kappa(\mu(i)-x)}{rh} -\frac{ \sum_{j=0}^{M_\xi} d_j\gamma(i)xz_j}{rh}- \frac{u}{rl} +\frac{\Gamma}{r} +\sum_{j\ne i} \frac{q_{ij}}{r}
\eeq
\beq{coefficient2}
\ad a(x,i) = \frac{\sigma^2(i)}{2rh^2} +  \frac{\kappa(\mu(i)-x)}{rh},
\eeq
\beq{coefficient3}
\ad b(x,i)  =\frac{\sigma^2(i)}{2rh^2}.
 \eeq
Note that equation (\ref{discrete}) is equivalent to $ V(s,x,y,i) = {\cal F}_\xi(V)(s,x,y,i;h,k,l)$, it suffices to show the operator ${\cal F}_\xi$ has a fixed point.
Using the fact that the difference of sups is less than the sup of differences.  If we have two bounded functions $V, W $  defined on $ {\cal D}_R$, it is clear that 
\bea
\ad |{\cal F}_\xi(V)(s,x,y,i;h,k,l) - {\cal F}_\xi(W)(s,x,y,i;h,k,l)| \\
\ad \leq  \bigg |  \sup_{u \in U} \bigg[  \bigg( a(x,i;u) + b(s,x,i;u) - c(s,x,y,i;u) + \frac{1}{rk}  +  \sum_{j=0}^{N_\xi}\frac{c_j}{r}  + \sum_{j\ne i} \frac{q_{ij}}{r}(u)\\
\ad - \frac{u}{rl}-\frac{ \sum_{j=0}^{M_\xi} d_j\gamma(i)xz_j}{rh}\bigg)\sup_{{\cal D}_R}|V-W|\bigg] \bigg|\\
\ad \leq \bigg|\sum_{j=0}^{N_\xi}\frac{c_j}{r}-\frac{\Gamma}{r}\bigg|\sup_{{\cal D}_R}|V-W|.
\eea
Therefore, for $\xi \in(0,1) $ small enough so that $\disp \bigg|\sum_{j=0}^{N_\xi}\frac{c_j}{r}-\frac{\Gamma}{r}\bigg| <1$, the map ${\cal F}_\xi$ is a contraction on the space of bounded functions on $ {\cal D}_R$, using the Banach's Fixed Point Theorem we conclude the proof of the lemma.
\endproof

\subsection{Proof of \thmref{convergence}}
Define
\beq{sup} \barray
V^{*}(s,x,y,i)\ad=\limsup_{\theta \to s,\eta \to
x,\zeta\to y, k\downarrow 0, h \downarrow 0, l\downarrow 0}
V_{k,h, l}(\theta,\eta,\zeta,i) \,\, \hbox{and}  \\
V_{*}(s,x,y,i)\ad=\liminf_{\theta \to s,\eta \to
x,\zeta\to y, k\downarrow 0, h \downarrow 0, l\downarrow 0}
V_{k,h, l}(\theta,\xi,\zeta,i) \,\, .
\earray\eeq
We claim that $ V^{*}$ and $V_{*}$
are  sub- and supersolutions of (\ref{isaac1}), respectively.

To prove this claim, we only
consider the  case for $V^*$.
The argument for that of $V_{*}$ is similar.
 For each $i\in\M$, we want to show
\[
\disp
\frac{\partial \Phi}{\partial s}(s_0,x_0,y_0) \leq  H(s_0,x_0,y_0,i,V^*,\frac{\partial \Phi}{\partial x}, \frac{\partial \Phi}{\partial y} , \frac{\partial^2 \Phi}{\partial x^2}),
\]  for any test function $\Phi \in {\cal C}^{1,2,1}([0,T]\times \rr^+\times [0,K]) $
such that $(s_0,x_0,y_0, i)$ is a strictly local
 maximum of $V^* (s,x,y,i) -\Phi(s,x,y) $.
 Without loss of generality, we may assume
 that $V^*(s_0,x_0,y_0,i) = \Phi(s_0,x_0,y_0)$
 and because of the stability of our scheme we
 can also assume that
 $\Phi \geq  \sup_{k,h,l}\|V_{k,h,l}\|$ outside of the ball
 $ B((s_0,x_0,y_0),r) $
 where $r>0$ is such that
\[ \disp
V^*(s,x,y,i)- \Phi(s,x,y)\leq 0=V^*(s_0,x_0,y_0,i)-
\Phi(s_0,x_0,y_0) \ \hbox{ in } \  B((s_0,x_0,y_0),r).
\]
This implies that there exist sequences $k_n>0$,
$h_n >0$, $l_n>0$  and $(\theta_n, \eta_n,\zeta_n) \in[0,T]\times \rr^+\times[0,K]$
such that as $n\rightarrow \infty$ we have
\beq{limit3} \barray
\ad k_n \rightarrow 0, \,\,\,h_n\rightarrow 0,
\,\,\, l_n\to 0,\,\,\, \theta \to s_0,\,\,\,  \eta_n\rightarrow x_0,\, \,\,\zeta_n\rightarrow y_0,\,\,\,
V_{k_n,h_n,l_n}(\theta_n, \eta_n,\zeta_n,i)\rightarrow V^*(s_0,x_0,y_0,i),\, \nonumber \\
\ad \hbox{and}\,\,\, (\theta_n, \eta_n,\zeta_n) \,
\hbox{ is  a global maximum of } V_{k_n,h_n, l_n}-\Phi. \nonumber
\earray\eeq
Denote $\epsilon_n=V_{k_n,h_n, l_n}
(\theta_n,\eta_n,\zeta_n,i)-\Phi(\theta_n,\eta_n,\zeta_n)$. Obviously
$\epsilon_n \rightarrow 0$ and
\beq{limit4}V_{k_n,h_n,l_n}(s,x,y,i)\leq \Phi(s,x,y)
+ \epsilon_n  \ \hbox{ for all } \ (s,x,y) \in [0,T]\times  \rr^+\times [0,K].
\eeq
We know that for all $\xi \in(0,1)$,
  \[
S(\xi,k_n,h_n,l_n,\theta_n,\eta_n,\zeta_n,i,V_{k_n,h_n, l_n}(\theta_n,\eta_n,\zeta_n,i),V_{k_n,h_n,l_n})=0
.\] The monotonicity of $S$ and (\ref{limit4}) implies
\beq{ineq3}  \barray \ad
S(\xi,k_n,h_n,l_n,\theta_n,\eta_n,\zeta_n,i,\Phi(\theta_n,\eta_n,\zeta_n) + \epsilon_n,\Phi + \epsilon_n)   \\
\aad  \ \le S(\xi,k_n,h_n,l_n,\theta_n,\eta_n,\zeta_n,i,V_{k_n,h_n, l_n}(\theta_n,\eta_n,\zeta_n,i),V_{k_n,h_n,l_n})=0.
\earray
\eeq
Therefore,
\[
\disp
\lim_{\xi\downarrow 0}\lim_{n\to \infty}S(\xi,k_n,h_n,l_n,\theta_n,\eta_n,\zeta_n,i,\Phi(\theta_n,\eta_n,\zeta_n) + \epsilon_n,\Phi + \epsilon_n)  \leq 0
,\]
so
\beq{ineq3000}\barray \ad
\frac{\partial \Phi}{\partial s}(s_0,x_0,y_0) \leq  H(s_0,x_0,y_0,i,V^*,\frac{\partial \Phi}{\partial x}, \frac{\partial \Phi}{\partial y} , \frac{\partial^2 \Phi}{\partial x^2}). \nonumber
\earray \eeq
This proves that $V^*$ is a viscosity subsolution and,
similarly we can prove that $V_{*} $ is  a viscosity
supersolution. Thus, using the uniqueness of
the viscosity solution, we see that $V=V^*=V_{*}$. Therefore,
we
conclude that the sequence $(V_{h,k,l})_{h,k,l}$
converges  locally uniformly  to $V$ as desired.
\endproof

\end{document}